\documentclass[11pt]{amsart}
\usepackage{amsmath, amsthm, amssymb}

\usepackage[a4paper,width=150mm,top=25mm,bottom=25mm,centering]{geometry}

\usepackage{hyperref} 
\usepackage{enumitem}
\usepackage{graphicx} 
\usepackage{xcolor}
\usepackage{comment}
\usepackage{url}
\usepackage{caption} 

\newtheorem{theorem}{Theorem}[section]
\newtheorem{conjecture}[theorem]{Conjecture}

\newtheorem{corollary}[theorem]{Corollary}

\theoremstyle{definition}

\title[Non-fibered twisted torus knots]{Infinite families of non-fibered twisted torus knots}

\author{Adnan}
\address{Department of Mathematics, Kangwon National University, Gangwon 24341,
 Republic of Korea}
\email{adnanshahab35@kangwon.ac.kr}

\author{Kyungbae Park}
\address{Department of Mathematics, Kangwon National University, Gangwon 24341,
 Republic of Korea}
\email{kyungbaepark@kangwon.ac.kr}

\subjclass[2020]{57K10, 57K14}
\keywords{twisted torus knots, fibered knots, Alexander polynomial}

\begin{document}
\begin{abstract}
    We present explicit infinite families of twisted torus knots that are not fibered. Our approach relies on an explicit formula for the Alexander polynomial derived in our previous work. We show that the leading coefficients of the Alexander polynomials of twisted torus knots can take arbitrary integer values, which immediately yields infinitely many examples of non-fibered twisted torus knots.
\end{abstract}

\maketitle

\section{Introduction}
\begin{figure}[t]
    \centering   
    \includegraphics[angle=90, width=0.6\textwidth]{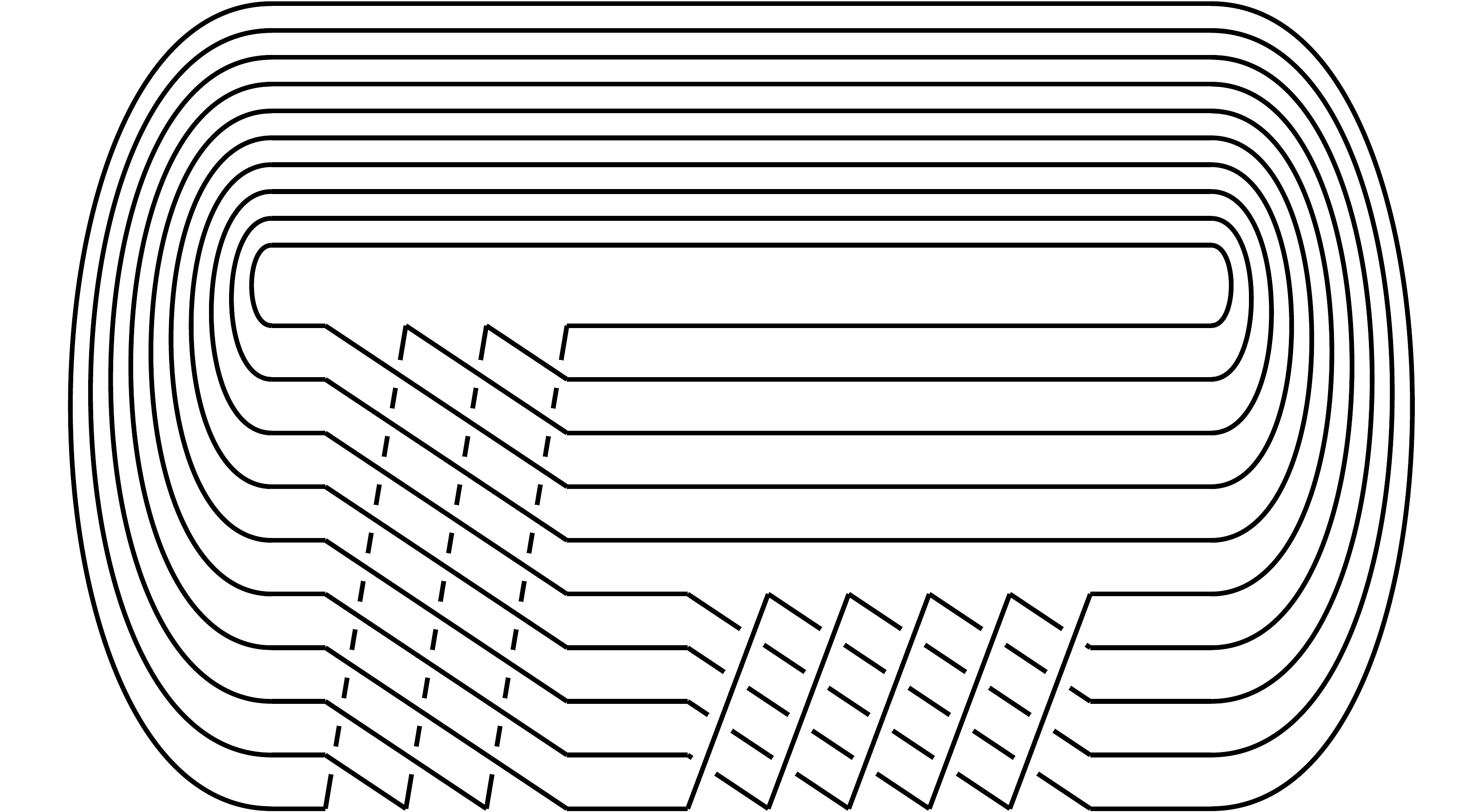}
    \caption{The twisted torus knot $T(10,3;5,-1)$, which is not fibered.}
    \label{Figure-TTK}
\end{figure}

For relatively prime integers $p$ and $q$, and integers $0 < r \leq |p|+|q|$ and $s\in\mathbb{Z}$, a \emph{twisted torus knot} $T(p,q;r,s)$ is a generalization of torus knots obtained by applying $s$ full twists on $r$ adjacent strands of the standard diagram of the $(p,q)$-torus knot on the flat torus. Figure~\ref{Figure-TTK} illustrates the knot $T(10,3;5,-1)$. 

Twisted torus knots, introduced by Dean in his thesis \cite{Dean-1996}, have been studied from many geometric and topological perspectives. These include hyperbolic volume \cite{CDW-1999, CKP-2004, CFKNP-2011}, Seifert fibered and graph-manifold surgeries \cite{Dean-2003, Kang-2016}, group-theoretic properties of knot groups \cite{Motegi-Teragaito-2022, Himeno-Teragaito-2023}, Heegaard splittings \cite{Moriah-Sedgwick-2009}, $L$-space surgeries and related questions arising from the $L$-space conjecture \cite{Vafaee-2015, CGHV-2016, Ichihara-Yuki-2018, Tran-2020}, tunnel number \cite{JH_Lee-2011}, and connections with Lorenz knots (or T-links) \cite{Birman-Kofman-2009}. In addition, the problem of classifying twisted torus knots as torus, satellite, or hyperbolic has been studied extensively \cite{Lee-2014, Lee-2015, Lee-2018, Lee-2019, Lee-2021, Lee-Paiva-2022}.

In this paper, we study the fiberedness of twisted torus knots. Recall that a knot $K \subset S^3$ is \emph{fibered} if its exterior $S^3 \setminus \nu(K)$ fibers over the circle $S^1$ with fiber a Seifert surface for $K$ \cite{Stallings-1978}. In particular, torus knots form a classical family of fibered knots. Fibered knots play a central role in low-dimensional topology, as their exteriors admit particularly simple topological structures encoded by the monodromy of the fiber. They also play an important role in the study of contact structures, knot Floer homology, and the interplay between algebraic and geometric knot invariants \cite{Etnyre-2006, Ni-2007, Ozsvath-Szabo-2005}.

Throughout this paper, we assume $p>q>0$, since $T(p,q;r,s)=T(q, p;r,s)$ and $T(p,q;r,s)$ is the mirror of $T(p,-q;r,-s)$. Since fiberedness is preserved under taking mirrors, this assumption causes no loss of generality. If $r \leq p$ and $s > 0$, then $T(p,q;r,s)$ is a positive braid knot, and hence fibered by Stallings’ theorem \cite{Stallings-1978}. On the other hand, when $s < 0$, the situation is more subtle. In this case, some twisted torus knots are known to be torus knots \cite{Lee-2021}, and some are even the unknot \cite{Lee-2014}; in particular, such knots are fibered. In \cite{Doleshal-2013}, Doleshal showed that if $q|s| < p$ and $r < q$, then $T(p,q;r,s)$ is a positive braid knot, and hence fibered. 
However, Doleshal’s condition is not necessary for a twisted torus knot to be fibered. For example, Lee’s result \cite[Theorem~1.1]{Lee-2015} and the work of Lee and Paiva \cite[Theorem~1.1]{Lee-Paiva-2022} show that certain families of twisted torus knots with $q|s| \ge p$ or $r > q$ are in fact torus knots, and hence fibered.

In his paper, Doleshal pointed out that not all twisted torus knots are fibered; in particular,  
$T(4,3;2,-2)$ is not fibered, since it is isotopic to the five-crossing twist knot. In this paper, we exhibit explicit two-parameter families of non-fibered twisted torus knots. 

Our method relies on an explicit formula for the Alexander polynomial of twisted torus knots with $r\leq p$, obtained in our previous work \cite{Park-Adnan-2025-01}. It is well known that every fibered knot has a monic Alexander polynomial. In particular, we show that the leading coefficients of the Alexander polynomials of twisted torus knots can attain arbitrary integer values, and hence establish the existence of infinitely many non-fibered examples. Since the Alexander polynomial is defined only up to multiplication by \(\pm t^k\), we interpret the leading coefficient as the absolute value of the coefficient of the highest-degree term, and the degree as the difference between the highest and lowest exponents.

\begin{theorem}\label{thm:main-1}
    Let $r > 0$ and $s < -1$. Then the Alexander polynomial of the twisted torus knot 
    \[
        T(r|s|,\, r|s|-1;\, r,\, s)
    \]
    has leading coefficient $r$ and degree $r|s|(r|s| - r - 2) + 2$.
\end{theorem}

In particular, if $r > 1$, then the Alexander polynomials of these knots are not monic, and hence the knots are not fibered. Moreover, these knots are pairwise distinct, as they are distinguished by the degrees and leading coefficients of their Alexander polynomials. This yields the following corollary.

\begin{corollary}
The family of twisted torus knots
\[
    \{\,T(r|s|,\, r|s|-1;\, r,\, s) \mid r>1,\ s<-1\,\}
\]
consists of pairwise distinct non-fibered knots.
\end{corollary}

When the last parameter $s$ of a twisted torus knot is $-1$, the above family has monic Alexander polynomials. Nevertheless, we can still construct analogous one-parameter families in this case.

\begin{theorem}\label{thm:main-2}
For each integer $n > 0$, the Alexander polynomial of the twisted torus knot
\[
    T(6n-1,\, 2n;\, 3n,\,-1) 
\]
has leading coefficient $n$ and degree $(3n-2)(n-1)$. 
\end{theorem}

We therefore obtain the following corollary.

\begin{corollary}
The family of twisted torus knots
\[
    \{\,T(6n - 1,\, 2n;\, 3n,\,-1) \mid n > 1\,\}
\]
consists of pairwise distinct non-fibered knots.
\end{corollary}

In fact, similar computations of the Alexander polynomial yield the following additional families in the cases $s=-1$ and $s=-2$.
\begin{theorem}\label{thm:main-3}
    For each integer $n > 1$, the Alexander polynomials of the following twisted torus knots are not monic. In particular, these knots are not fibered:
    \begin{enumerate}
        \item $T(6n-2,\, 2n-1;\, 3n-1,\,-1)$
        \item  $T(10n-6,\,2n-1;\,5n-3,\,-1)$
        \item  $T(10n+1,\,2n;\,5n,\,-1)$
        \item  $T(4n+4,\,2n+1;\, 2n+2,\,-2)$
        \item  $T(6n+1,\, 2n ;\, 3n,\,-2)$
        \item  $T(6n+2,\, 2n+1;\, 3n+1,\,-2)$
        \item  $T(18n+12,\,6n+5;\, 9n+6,\,-2)$
        \item  $T(18n+18,\,6n+7;\, 9n+9,\,-2)$
    \end{enumerate}
\end{theorem}

Since every $L$-space knot is fibered, our examples also provide new explicit families of twisted torus knots with $s < 0$ that are not $L$-space knots. For previous examples of twisted torus knots that are not $L$-space knots, see \cite[Theorem 3 and 4]{Morton-2006}, \cite[Theorem 1.1]{Vafaee-2015} and \cite[Theorem 3]{Park-Adnan-2025-01}.

Although we have found many examples of non-fibered twisted torus knots, a complete classification appears to be considerably more difficult. Beyond the Alexander polynomial, there are other ways to detect fiberedness. For instance, Stallings \cite{Stallings-1978} proved that a knot is fibered if and only if the commutator subgroup of its knot group is finitely generated (equivalently, free of finite rank). In our previous work \cite{Park-Adnan-2025-01}, we obtained a presentation of the knot group of a twisted torus knot, although it seems difficult to apply Stallings' criterion directly in this setting.

A further refinement of the group-theoretic approach is provided by twisted Alexander polynomials, which are defined using linear representations of the knot group and often yield stronger obstructions to fiberedness than the ordinary Alexander polynomial. In particular, Cha showed that twisted Alexander invariants of a fibered knot satisfy strong algebraic constraints \cite{Cha-2003}, and Goda, Kitano, and Morifuji proved that, for suitable representations, the associated twisted Alexander polynomial gives a necessary condition for a knot in $S^3$ to be fibered \cite{Goda-Kitano-Morifuji-2005}. For a general overview of this theory and its applications, see \cite{Friedl-Vidussi-2011}.

Another important approach is through knot Floer homology \cite{Ozsvath-Szabo-2004, Rasmussen-2003}, which provides a powerful tool for detecting fiberedness. A fundamental theorem of Ghiggini \cite{Ghiggini-2008} and Ni \cite{Ni-2007} states that a knot $K\subset S^3$ is fibered if and only if
\[
\operatorname{rank}\,\widehat{HFK}(K,g(K))=1,
\]
where $g(K)$ denotes the Seifert genus of $K$. Despite the strength of these criteria, applying them systematically to twisted torus knots appears to be difficult in general.

On the other hand, it is well known that monicity of the Alexander polynomial is not, in general, sufficient to guarantee fiberedness. Notable examples include the Kinoshita--Terasaka knot and the Conway knot, which are mutants of one another. Both have trivial (and hence monic) Alexander polynomials, but neither is fibered \cite{Goda-Kitano-Morifuji-2005, Ozsvath-Szabo-2005}. Nevertheless, our computations suggest that twisted torus knots may behave differently. Among the $2{,}152$ twisted torus knots with $s<0$ whose standard braid diagrams have at most $100$ crossings, we found, using knot Floer homology computations in \texttt{SnapPy} \cite{SnapPy}, that $49$ are non-fibered. We further found that the Alexander polynomials of all these knots are non-monic. Although the evidence is admittedly limited, these computations suggest that fiberedness for twisted torus knots may be governed by a simpler principle. This leads us to formulate the following conjecture, which is known to hold for alternating knots by a classical theorem of Murasugi \cite{Murasugi-1958}.

\begin{conjecture}
A twisted torus knot is fibered if and only if its Alexander polynomial is monic.
\end{conjecture}

\subsection*{Acknowledgments}
The authors were supported by the National Research Foundation of Korea (NRF) grants funded by the Korean government (RS-2025-24523511 and RS-2025-25415913).

\section{Alexander polynomial of twisted torus knots}
In this section, we recall, for convenience, the formula for the Alexander polynomial of twisted torus knots from \cite{Park-Adnan-2025-01}. We refer the reader to that paper for full details and proofs.  

For relatively prime integers $p>q>0$, and integer $r$ with $1<r<p$, and $s\in\mathbb{Z}$, we first introduce some notation related to the modular arithmetic of $p$, $q$, and $r$. Let $[x]$ denote the residue of $x$ modulo $p$; that is, the unique integer satisfying $0 \leq [x] < p$ and $[x]\equiv x$ modulo $p$. 

We first define
\[
Q:=\{[jq^{-1}] \mid j=0,1,\dots,r-1\},
\]
where $q^{-1}$ denotes the multiplicative inverse of $q$ modulo $p$. We order the elements of $Q$ increasingly and write
\[
Q=\{Q_0<Q_1<\cdots<Q_{r-1}\}.
\]

Next, let
\[
R:=\{[-jq^{-1}] \mid j=1,2,\dots,q\}.
\]
For each $i=1,\dots,r-1$, we define $k_i$ to be the number of elements of $R$ lying in the interval $[Q_{i-1},Q_i)$; that is,
\[
k_i:=\bigl|\{x\in R \mid Q_{i-1}\le x<Q_i\}\bigr|.
\]
 
Define
\[
Q':=[rq^{-1}].
\]
Since $p$ and $q$ are relatively prime, we have $Q'\notin Q$. Hence there exists an integer $m$ with $0\le m\le r-1$ such that
\[
Q_0<\cdots<Q_m<Q'<Q_{m+1}<\cdots<Q_{r-1}.
\]
We then define
\[
k':=\bigl|\{x\in R \mid Q_m\le x<Q'\}\bigr|.
\]

For each $i=1,\dots,r-1$, we also define
\[
\overline{k}_i:=k_1+\cdots+k_i,
\]
and set
\[
\overline{k}':=k_1+\cdots+k_m+k'.
\]

With the above notation, the Alexander polynomial of a twisted torus knot is given by the following formula \cite{Park-Adnan-2025-01}.
\begin{equation}\label{eq:Alexander_poly}
       \Delta_{T(p,q;r,s)}(t)=\frac{1-t}{(1-t^p)(1-t^q)(1-t^{r})}\left(\tilde{X}(t)Y(t)-X(t)\tilde{Y}(t)\right), 
\end{equation}
where
\begin{align*}
   X(t)
   &=1-(1-t^{rs})\sum_{i=1}^{r-1}t^{\overline{k}_ip+(i-1)rs}-t^{pq+(r-1)rs},\\
   \widetilde{X}(t)&=\begin{cases}
        1-t^{\overline{k'}p},&\text{if } m=0,\\
        1-(1-t^{rs})\displaystyle\sum_{i=1}^{m}t^{\overline{k}_i p+(i-1)rs}-t^{\overline{k'}p+mrs},&\text{otherwise,}    
   \end{cases}\\
   Y(t)&=1-(1-t^{rs})\sum_{i=1}^{r-1}t^{Q_iq+(i-1)rs}-t^{pq+(r-1)rs},\\
   \widetilde{Y}(t)&=\begin{cases}
       1-t^{Q'q},&\text{if } m = 0,\\
        1-(1-t^{rs})\displaystyle\sum_{i=1}^{m}t^{Q_iq+(i-1)rs}-t^{Q'q+mrs},&\text{otherwise.}       \end{cases}
\end{align*}

\section{Leading coefficients and degrees of certain twisted torus knots}
Using the notation from the previous section, we determine the leading coefficients and degrees of the Alexander polynomials for the families appearing in Theorems \ref{thm:main-1} and \ref{thm:main-2}. These computations establish those theorems. Theorem \ref{thm:main-3} can be proved similarly, and we leave the proof to the reader.

\subsection{Proof of Theorem \ref{thm:main-1}}
For the family $T(r|s|,r|s|-1;r,s)$ with $s<-1$, we have
\[
p=r|s| \qquad\text{and}\qquad q=r|s|-1\equiv -1 \pmod p.
\]
Hence
\[
[q^{-1}]=r|s|-1.
\]
Therefore
\[
Q=\{[j(r|s|-1)] \mid j=0,1,\dots,r-1\},
\]
and, after arranging the elements increasingly, we obtain
\[
Q_i=
\begin{cases}
0 & \text{if } i=0,\\
r|s|-r+i & \text{if } 1\le i\le r-1.
\end{cases}
\]

Moreover,
\[
R=\{[-j(r|s|-1)] \mid j=1,2,\dots,r|s|-1\}
 =\{1,2,\dots,r|s|-1\}.
\]
It follows that
\[
k_i=
\begin{cases}
r|s|-r & \text{if } i=1,\\
1 & \text{if } 2\le i\le r-1.
\end{cases}
\]

Since
\[
Q'=[r(r|s|-1)]=r|s|-r,
\]
we have
\[
Q_0<Q'<Q_1.
\]
Hence \(m=0\), so that \(Q_m=Q_0\), and
\[
k'=r|s|-r-1.
\]
Therefore,
\[
\overline{k}_i=r|s|-r-1+i \qquad \text{for } 1\le i\le r-1,
\]
and
\[
\overline{k}'=r|s|-r-1.
\]

Hence, the polynomials \(X(t)\), \(\widetilde{X}(t)\), \(Y(t)\), and \(\widetilde{Y}(t)\) defined in the previous section are given by
\begin{align*}
X(t)
&=1-r t^{(r|s|-r)r|s|}+(r-1)t^{(r|s|-r-1)r|s|},\\
\widetilde{X}(t)
&=1-t^{(r|s|-r-1)r|s|},\\
Y(t)
&=1-(1-t^{rs})\sum_{i=1}^{r-1} t^{r|s|(r|s|-r)+r-i}-t^{r|s|(r|s|-r)},\\
\widetilde{Y}(t)
&=1-t^{(r|s|-r)(r|s|-1)}.
\end{align*}

We obtain
\begin{align*}
\widetilde{X}(t)Y(t)-X(t)\widetilde{Y}(t)
={}&-r t^{r|s|(2r|s|-2r-1)+r} \\
&+(t^{r|s|(r|s|-r-1)}-1)(1-t^{rs})\sum_{i=1}^{r-1} t^{r|s|(r|s|-r)+r-i} \\
&+t^{r|s|(2r|s|-2r-1)} +(r-1)t^{2r|s|(r|s|-r-1)+r} \\
&+(r-1)t^{r|s|(r|s|-r)} + t^{(r|s|-r)(r|s|-1)} - r t^{r|s|(r|s|-r-1)}.
\end{align*}

To determine the leading term of \(\Delta_{T(r|s|,\,r|s|-1;\,r,\,s)}(t)\), we expand the rational factors formally as geometric series:
\[
\Delta_{T(r|s|,\,r|s|-1;\,r,\,s)}(t)
=
\frac{1-t}{(1-t^{r|s|})(1-t^{r|s|-1})(1-t^r)}
\bigl(\widetilde{X}(t)Y(t)-X(t)\widetilde{Y}(t)\bigr).
\]
Thus, as a formal power series,
\begin{align*}
\Delta_{T(r|s|,\,r|s|-1;\,r,\,s)}(t)
={}&
\bigl(\widetilde{X}(t)Y(t)-X(t)\widetilde{Y}(t)\bigr)(1-t) \\
&\cdot (1+t^{r|s|}+t^{2r|s|}+\cdots)
(1+t^{r|s|-1}+t^{2(r|s|-1)}+\cdots)
(1+t^r+t^{2r}+\cdots).
\end{align*}
We see that the lowest-degree term is
\[
-r\,t^{(r|s|-r-1)r|s|}.
\]
Since the Alexander polynomial is symmetric, the coefficient of the highest-degree term is also \(-r\). Hence the leading coefficient of $\Delta_{T(r|s|,\,r|s|-1;\,r,\,s)}(t)$ is \(r\).

The lowest exponent is determined by the lowest-degree term in the formal power series expansion at \(t=0\). The highest exponent, however, is more conveniently obtained by comparing the degrees of the numerator and denominator in the rational expression for the Alexander polynomial. Hence, the highest and lowest exponents of $\Delta_{T(r|s|,\,r|s|-1;\,r,\,s)}(t)$
are
\[
2(r|s|)^2-2r^2|s|-3r|s|+2
\quad\text{and}\quad
r|s|(r|s|-r-1),
\]
respectively. Therefore, the degree of $\Delta_{T(r|s|,\,r|s|-1;\,r,\,s)}(t)$ is $r|s|(r|s|-r-2)+2$.

\subsection{Proof of Theorem \ref{thm:main-2}}
For the twisted torus knot \(T(6n-1,\,2n;\,3n,\,-1)\), we have
\[
[q^{-1}]=[(2n)^{-1}]=3.
\]
Hence
\[
Q=\{[3j] \mid j=0,1,2,\dots,3n-1\}=\{Q_0<Q_1<\cdots<Q_{3n-1}\}.
\]
A straightforward computation shows that
\[
Q_i=
\begin{cases}
0 & \text{if } i=0,\\[2mm]
i+\dfrac{i-1}{2} & \text{if } 1\le i\le 2n \text{ and } i \text{ is odd},\\[2mm]
i+\dfrac{i}{2} & \text{if } 1\le i\le 2n \text{ and } i \text{ is even},\\[2mm]
3(i-n) & \text{if } 2n+1\le i\le 3n-1.
\end{cases}
\]

Similarly,
\[
R=\{[-3j] \mid j=1,2,\dots,2n\}=\{R_1<R_2<\cdots<R_{2n}\},
\]
where
\[
R_i=
\begin{cases}
3i-1 & \text{if } 1\le i\le 2n-1,\\
6n-2 & \text{if } i=2n.
\end{cases}
\]
It follows that
\[
k_i=
\begin{cases}
0 & \text{if } 1\le i\le 2n-1 \text{ and } i \text{ is odd},\\
1 & \text{if } 1\le i\le 2n-1 \text{ and } i \text{ is even},\\
1 & \text{if } 2n\le i\le 3n-1.
\end{cases}
\]

Since $Q'=[3(3n)]=3n+1$, we have
\[
Q_{2n}<Q'<Q_{2n+1}.
\]
Hence \(m=2n\), so that \(Q_m=Q_{2n}\), and $k'=0$. Therefore,
\[
\overline{k}_i=
\begin{cases}
\left\lfloor \dfrac{i}{2} \right\rfloor & \text{if } 1\le i\le 2n,\\[2mm]
i-n & \text{if } 2n+1\le i\le 3n-1,
\end{cases}
\]
and
\[
\overline{k}'=n.
\]
Here \(\lfloor x\rfloor\) denotes the floor function.

Thus, we have
\begin{align*}
X(t)
&= t^{2n-1}-\sum_{i=0}^{n-1} t^{2n+i}
 +(1-t)\sum_{i=0}^{n-3} t^{5n-2+i(3n-1)}
 +t^{-n}+\sum_{i=0}^{n-1} t^{-3n-i} \\
&\qquad -t^{3n^2-2n+1}-t^{3n^2+n},\\[1mm]
\widetilde{X}(t)
&= \sum_{i=0}^{n-1} t^{-3n-i}-\sum_{i=0}^{n-1} t^{2n+i},\\[1mm]
Y(t)
&= nt^{-n}+(n+1)-nt^{2n}-(n-1)t^{3n}-t^{3n^2}-t^{3n^2+n},\\[1mm]
\widetilde{Y}(t)
&= nt^{-n}+(n+1)-(n+1)t^{2n}-nt^{3n}.
\end{align*}
Therefore,
\begin{small}
\begin{align*}
\widetilde{X}(t)Y(t)-X(t)\widetilde{Y}(t)
={}&-n t^{3n^2+4n}+n t^{3n^2-3n+1}+n t^{3n^2}-n t^{3n^2+n+1}\\
&+\bigl(t^{3n^2+n}+t^{3n^2}-t^{2n}-t^{3n}\bigr)\sum_{i=0}^{n-1} t^{2n+i}\\
&+(n+1)\bigl(t^{3n^2+n}+t^{3n^2-2n+1}+t^n-t^{3n^2+1}-t^{3n^2+3n}+t^{4n-1}-t^{2n-1}-t^{-n}\bigr)\\
&+\bigl((n+1)t^{2n}+n t^{3n}-(n+1)-n t^{-n}\bigr)(1-t)\sum_{i=0}^{n-3} t^{5n-2+i(3n-1)}\\
&+\bigl(t^{2n}+t^{3n}-t^{3n^2}-t^{3n^2+n}\bigr)\sum_{i=0}^{n-1} t^{-3n-i}\\
&+n t^{5n-1}+n t^{2n}-n t^{n-1}-n t^{-2n}.
\end{align*}
\end{small}

From the formal power series expansion of Equation \eqref{eq:Alexander_poly}, we see that the lowest-degree term of the Alexander polynomial of $T(6n-1,\,2n;\,3n,\,-1)$ is $-n\,t^{-2n}$. By symmetry of the Alexander polynomial, this implies that the leading coefficient is \(n\). Moreover, the highest and lowest exponents of
\[
\Delta_{T(6n-1,\,2n;\,3n,\,-1)}(t)
\]
are \(3n^2-7n+2\) and \(-2n\), respectively. Therefore, its degree is $(3n-2)(n-1)$.

\clearpage

\bibliographystyle{alpha}

\bibliography{references-ttk}

\newcommand{\etalchar}[1]{$^{#1}$}
\begin{thebibliography}{CGHV16}

\bibitem[BK09]{Birman-Kofman-2009}
Joan Birman and Ilya Kofman.
\newblock A new twist on {L}orenz links.
\newblock {\em J. Topol.}, 2(2):227--248, 2009.

\bibitem[CDGW]{SnapPy}
Marc Culler, Nathan~M. Dunfield, Matthias Goerner, and Jeffrey~R. Weeks.
\newblock Snap{P}y, a computer program for studying the geometry and topology
  of $3$-manifolds.
\newblock Available at \url{http://snappy.computop.org}.

\bibitem[CDW99]{CDW-1999}
Patrick~J. Callahan, John~C. Dean, and Jeffrey~R. Weeks.
\newblock The simplest hyperbolic knots.
\newblock {\em J. Knot Theory Ramifications}, 8(3):279--297, 1999.

\bibitem[CFK{\etalchar{+}}11]{CFKNP-2011}
Abhijit Champanerkar, David Futer, Ilya Kofman, Walter Neumann, and Jessica~S.
  Purcell.
\newblock Volume bounds for generalized twisted torus links.
\newblock {\em Math. Res. Lett.}, 18(6):1097--1120, 2011.

\bibitem[CGHV16]{CGHV-2016}
Katherine Christianson, Justin Goluboff, Linus Hamann, and Srikar Varadaraj.
\newblock Non-left-orderable surgeries on twisted torus knots.
\newblock {\em Proc. Amer. Math. Soc.}, 144(6):2683--2696, 2016.

\bibitem[Cha03]{Cha-2003}
Jae~Choon Cha.
\newblock Fibred knots and twisted {A}lexander invariants.
\newblock {\em Trans. Amer. Math. Soc.}, 355(10):4187--4200, 2003.

\bibitem[CKP04]{CKP-2004}
Abhijit Champanerkar, Ilya Kofman, and Eric Patterson.
\newblock The next simplest hyperbolic knots.
\newblock {\em J. Knot Theory Ramifications}, 13(7):965--987, 2004.

\bibitem[Dea96]{Dean-1996}
John~Charles Dean.
\newblock {\em Hyperbolic knots with small {S}eifert-fibered {D}ehn surgeries}.
\newblock ProQuest LLC, Ann Arbor, MI, 1996.
\newblock Thesis (Ph.D.)--The University of Texas at Austin.

\bibitem[Dea03]{Dean-2003}
John~C. Dean.
\newblock Small {S}eifert-fibered {D}ehn surgery on hyperbolic knots.
\newblock {\em Algebr. Geom. Topol.}, 3:435--472, 2003.

\bibitem[Dol13]{Doleshal-2013}
Brandy~Guntel Doleshal.
\newblock Fibered and primitive/{S}eifert twisted torus knots.
\newblock {\em J. Knot Theory Ramifications}, 22(1):1250141, 14, 2013.

\bibitem[Etn06]{Etnyre-2006}
John~B. Etnyre.
\newblock Lectures on open book decompositions and contact structures.
\newblock In {\em Floer homology, gauge theory, and low-dimensional topology},
  volume~5 of {\em Clay Math. Proc.}, pages 103--141. Amer. Math. Soc.,
  Providence, RI, 2006.

\bibitem[FV11]{Friedl-Vidussi-2011}
Stefan Friedl and Stefano Vidussi.
\newblock A survey of twisted {A}lexander polynomials.
\newblock In {\em The mathematics of knots}, volume~1 of {\em Contrib. Math.
  Comput. Sci.}, pages 45--94. Springer, Heidelberg, 2011.

\bibitem[Ghi08]{Ghiggini-2008}
Paolo Ghiggini.
\newblock Knot {F}loer homology detects genus-one fibred knots.
\newblock {\em Amer. J. Math.}, 130(5):1151--1169, 2008.

\bibitem[GKM05]{Goda-Kitano-Morifuji-2005}
Hiroshi Goda, Teruaki Kitano, and Takayuki Morifuji.
\newblock Reidemeister torsion, twisted {A}lexander polynomial and fibered
  knots.
\newblock {\em Comment. Math. Helv.}, 80(1):51--61, 2005.

\bibitem[HT23]{Himeno-Teragaito-2023}
Keisuke Himeno and Masakazu Teragaito.
\newblock New families of hyperbolic twisted torus knots with generalized
  torsion.
\newblock {\em Bull. Korean Math. Soc.}, 60(1):203--223, 2023.

\bibitem[IT18]{Ichihara-Yuki-2018}
Kazuhiro Ichihara and Yuki Temma.
\newblock Non-left-orderable surgeries on negatively twisted torus knots.
\newblock {\em Proc. Japan Acad. Ser. A Math. Sci.}, 94(5):49--52, 2018.

\bibitem[Kan16]{Kang-2016}
Sungmo Kang.
\newblock Twisted torus knots with graph manifold {D}ehn surgeries.
\newblock {\em Bull. Korean Math. Soc.}, 53(1):273--301, 2016.

\bibitem[LdP22]{Lee-Paiva-2022}
Sangyop Lee and Thiago de~Paiva.
\newblock Torus knots obtained by negatively twisting torus knots.
\newblock {\em J. Knot Theory Ramifications}, 31(1):Paper No. 2150080, 19,
  2022.

\bibitem[Lee11]{JH_Lee-2011}
Jung~Hoon Lee.
\newblock Twisted torus knots {$T(p, q; 3, s)$} are tunnel number one.
\newblock {\em Journal of Knot Theory and Its Ramifications}, 20(06):807--811,
  2011.

\bibitem[Lee14]{Lee-2014}
Sangyop Lee.
\newblock Twisted torus knots that are unknotted.
\newblock {\em Int. Math. Res. Not. IMRN}, (18):4958--4996, 2014.

\bibitem[Lee15]{Lee-2015}
Sangyop Lee.
\newblock Torus knots obtained by twisting torus knots.
\newblock {\em Algebraic \& Geometric Topology}, 15(15):2817–2836, 2015.

\bibitem[Lee18]{Lee-2018}
Sangyop Lee.
\newblock Satellite knots obtained by twisting torus knots: hyperbolicity of
  twisted torus knots.
\newblock {\em Int. Math. Res. Not. IMRN}, (3):785--815, 2018.

\bibitem[Lee19]{Lee-2019}
Sangyop Lee.
\newblock Positively twisted torus knots which are torus knots.
\newblock {\em J. Knot Theory Ramifications}, 28(3):1950023, 13, 2019.

\bibitem[Lee21]{Lee-2021}
Sangyop Lee.
\newblock Twisted torus knots {$T(mn+m+1,mn+1,mn+m+2,-1)$} and
  {$T(n+1,n,2n-1,-1)$} are torus knots.
\newblock {\em J. Knot Theory Ramifications}, 30(3):Paper No. 2150016, 20,
  2021.

\bibitem[Mor06]{Morton-2006}
Hugh~R. Morton.
\newblock The {A}lexander polynomial of a torus knot with twists.
\newblock {\em J. Knot Theory Ramifications}, 15(8):1037--1047, 2006.

\bibitem[MS09]{Moriah-Sedgwick-2009}
Yoav Moriah and Eric Sedgwick.
\newblock Heegaard splittings of twisted torus knots.
\newblock {\em Topology Appl.}, 156(5):885--896, 2009.

\bibitem[MT22]{Motegi-Teragaito-2022}
Kimihiko Motegi and Masakazu Teragaito.
\newblock Generalized torsion for knots with arbitrarily high genus.
\newblock {\em Canad. Math. Bull.}, 65(4):867--881, 2022.

\bibitem[Mur58]{Murasugi-1958}
Kunio Murasugi.
\newblock On the {A}lexander polynomial of the alternating knot.
\newblock {\em Osaka Math. J.}, 10:181--189; errata, 11 (1959), 95, 1958.

\bibitem[Ni07]{Ni-2007}
Yi~Ni.
\newblock Knot {F}loer homology detects fibred knots.
\newblock {\em Invent. Math.}, 170(3):577--608, 2007.

\bibitem[OS04]{Ozsvath-Szabo-2004}
Peter Ozsv\'ath and Zolt\'an Szab\'o.
\newblock Holomorphic disks and knot invariants.
\newblock {\em Adv. Math.}, 186(1):58--116, 2004.

\bibitem[OS05]{Ozsvath-Szabo-2005}
Peter Ozsváth and Zoltán Szabó.
\newblock On knot {F}loer homology and lens space surgeries.
\newblock {\em Topology}, 44(6):1281--1300, 2005.

\bibitem[PA25]{Park-Adnan-2025-01}
Kyungbae Park and Adnan.
\newblock The {A}lexander polynomial of twisted torus knots.
\newblock {\em J. Knot Theory Ramifications}, 34(8):Paper No. 2550038, 23,
  2025.

\bibitem[Ras03]{Rasmussen-2003}
Jacob~Andrew Rasmussen.
\newblock {\em Floer homology and knot complements}.
\newblock ProQuest LLC, Ann Arbor, MI, 2003.
\newblock Thesis (Ph.D.)--Harvard University.

\bibitem[Sta78]{Stallings-1978}
John~R. Stallings.
\newblock Constructions of fibred knots and links.
\newblock In {\em Algebraic and geometric topology ({P}roc. {S}ympos. {P}ure
  {M}ath., {S}tanford {U}niv., {S}tanford, {C}alif., 1976), {P}art 2}, volume
  XXXII of {\em Proc. Sympos. Pure Math.}, pages 55--60. Amer. Math. Soc.,
  Providence, RI, 1978.

\bibitem[Tra20]{Tran-2020}
Anh~T. Tran.
\newblock Non-left-orderable surgeries on {L}-space twisted torus knots.
\newblock {\em Proc. Amer. Math. Soc.}, 148(1):447--456, 2020.

\bibitem[Vaf15]{Vafaee-2015}
Faramarz Vafaee.
\newblock On the knot {F}loer homology of twisted torus knots.
\newblock {\em Int. Math. Res. Not. IMRN}, (15):6516--6537, 2015.

\end{thebibliography}
\end{document}